\documentclass[12pt]{article}
\usepackage{hyperref}
\usepackage{cmap}                        
\usepackage[cp1251]{inputenc}            
\usepackage[english]{babel}
\usepackage[left=1in,right=1in,top=1in,bottom=1.5in]{geometry} 
\usepackage{amssymb,amsmath, amsthm, amscd,ifthen,amsfonts,latexsym,color,epsfig}
\usepackage{graphicx}
\usepackage[svgnames]{xcolor}
\usepackage{tikz}
\usetikzlibrary{patterns,decorations.text,positioning,arrows,shapes,decorations.markings}
\tikzstyle{vertex}=[circle,draw=black,fill=black,inner sep=0,minimum size=0.2cm,text=white,font=\footnotesize]

\date{}

\title{Tiling the plane with equilateral triangles}

\author{J\'anos Pach\thanks{EPFL, Lausanne and R\'enyi Institute, Budapest. Supported by Swiss National Science Foundation Grants 200020-162884 and 200021-165977. E-mail: {\tt pach@cims.nyu.edu}.} \and G\'abor Tardos\thanks{R\'enyi Institute, Budapest. Supported by the Cryptography ``Lend\"ulet'' project of the Hungarian Academy
of Sciences and by the National Research, Development and Innovation Office, NKFIH, projects K-116769 and SNN-117879.}}

\date{}

\begin{document}
\maketitle
\begin{abstract}
Let $\cal T$ be a tiling of the plane with equilateral triangles no two of which share a side. We prove that if the side lengths of the triangles are bounded from below by a positive constant, then $\cal T$ is periodic and it consists of translates of only at most three different triangles. As a corollary, we prove a theorem of Scherer and answer a question of Nandakumar. The same result has been obtained independently by Richter and Wirth.
\end{abstract}

\section{Introduction}
A collection of compact convex sets with nonempty interiors is said to {\em tile} a region $R$ if no two of them share an interior point and their union is equal to $R$.

In 1934, during his first visit to Trinity College, Cambridge, P. Erd\H os~\cite{Sch98} raised the following question: Is it possible to tile a unit square with finitely many smaller squares, no two of which are congruent? It has taken a few years before four promising students at Trinity College, R. L. Brooks, C. A. B. Smith, A. H. Stone, and W. T. Tutte~\cite{BSST40}, all of whom became prominent mathematicians, managed to answer Erd\H os's question in the affirmative. They discovered and explored an intimate relationship between tilings and flows in networks, and constructed several tilings with the required properties. In any such tiling $\cal T$ of a unit square $S$, there is a unique {\em smallest} square. This square cannot share a side with the boundary of $S$, otherwise one of its neighbors would be congruent to it. Scaling up $\cal T$ so that the side length of its smallest square becomes $1$, and tiling this unit square with a congruent copy of $\cal T$, we obtain a tiling ${\cal T}'\supset{\cal T}$ of a larger square with pairwise noncongruent squares. Repeating this procedure infinitely many times, we obtain the following.
\medskip

\noindent{\bf Corollary.} {\em The plane can be tiled by infinitely many pairwise noncongruent squares whose side lengths are bounded from below by a positive constant.}
\medskip

Tutte~\cite{Tu48} extended the technique developed in~\cite{BSST40} to tilings of an equilateral triangle with equilateral triangles. In particular, he proved that no such tiling exists with pairwise noncongruent triangles. However, this does not answer the corresponding question for tilings of the entire plane, which was raised by R. Nandakumar~\cite{Na14} (June 14, 2016).
\medskip

\noindent{\bf Problem 1.} (Nandakumar) {\em Is it possible to tile the plane with pairwise noncongruent equilateral triangles whose side lengths are bounded from below by a positive constant?}
\medskip

In fact, we do not even know whether there exists a tiling of the plane with pairwise noncongruent equilateral triangles if we replace the condition that they cannot be arbitrarily small by the weaker assumption that the tiling is {\em locally finite}, that is, every bounded region intersects only finitely many triangles.  Without assuming local finiteness, it is easy to find a tiling of the ``punctured'' plane and with a little more work of the whole plane; see \cite{Ri12}.

\medskip\noindent{\bf Proposition 2.} (Klaassen~\cite{Kl95}) {\em There exists a tiling of the plane minus a single point with equilateral triangles of pairwise distinct sizes.}
\medskip

\noindent{\it Proof}: Let $\alpha$ be the only real root of the polynomial
$x^3+x^2-1$. Let $P$ and $A$ be distinct points in the Euclidean plane. Let
$\phi$ be the
similarity transformation of the plane that has $P$ as its fixed point, rotates
the plane with angle $\pi/3$ and shrinks the distances with factor
$\alpha$. Let $T_0$ be the equilateral triangle with vertices $A$, $\phi(A)$ and
a third vertex in the half plane bounded by the line $A\phi(A)$ and containing
$P$.

The equilateral triangles $T_i=\phi^i(T_0)$, where
$i$ ranges over all integers, form a tiling of the plane minus the point $P$, and
the sizes of these triangles are pairwise distinct. $\Box$
\medskip

 \begin{figure}[ht]
\centering
\includegraphics[height=3.5in]{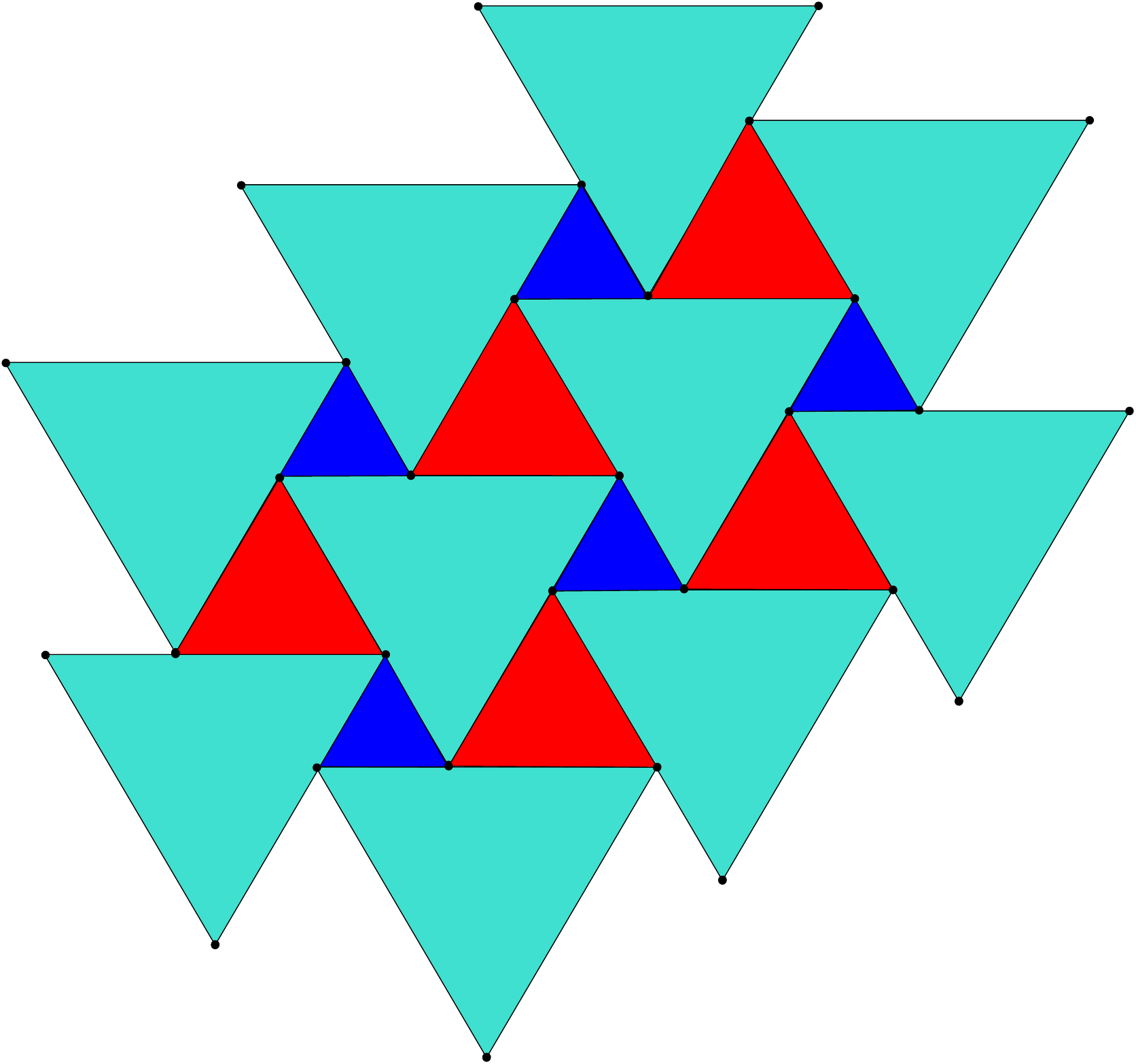}
\caption{A tiling with no two triangles sharing a side.}
\end{figure}

\smallskip

The aim of this note is to answer Nandakumar's above mentioned question (Problem 1).

\medskip\noindent{\bf Theorem 3.} {\em Let $\cal T$ be a tiling of the plane with equilateral triangles such that the side lengths of the triangles are bounded from below by a positive constant and no two triangles share a side. Then the triangles in $\cal T$ have at most three different side lengths, $a$, $b$ and $c$ with $a=b+c$ (where $c$ may be equal to $b$), and the tiling is periodic.}

\medskip\noindent{\bf Corollary 4.} {\em There is no tiling of the plane with pairwise noncongruent equilateral triangles whose side lengths are bounded from below by a positive constant.}
\medskip

The weaker statement that there is no tiling satisfying the conditions of Corollary 4, in which there is a {\em smallest} triangle, was proved by Scherer~\cite{Sch83}.  


Nandakumar~\cite{Na06, Na14}, who teaches undergraduates at a local college in Cochin, India, raised several other interesting questions about tilings. They have triggered a lot of research in geometry and topology; see~\cite{BBS10, BlZ14, Fr16, KaHA14, KuPT17a, KuPT17b, NaR12}.
\smallskip

After completing our note, it has come to our attention that Theorem 3 has also been proved independently and probably slightly earlier by Richter and Wirth~\cite{RiW17}. Their proof is somewhat different from ours, but in the last step we both use martingales. We are grateful to Richter and Wirth to their helpful remarks on our manuscript.

\section{Proof of Theorem 3}

Let us fix tiling $\cal T$ satisfying the requirements of Theorem 3. We will prove the statement of the theorem through a sequence of lemmas about $\cal T$. First, we agree about the terminology.

We call an edge of a triangle $T\in\cal T$ {\it subdivided} if some interior
point of this edge is a vertex of another triangle $T'\in\cal T$. Otherwise, it
is called {\it uncut}. Note that a vertex in the tiling is the interior point of at
most one edge and in that case it is the vertex of exactly three triangles.

We say that the edge $AB$ of a triangle {\it continues at $A$} if $A$ is the interior
point of an edge $e$ of another triangle of our tiling and $e$ also contains
some interior points of the edge $AB$. Clearly, at most one of the two edges
of a triangle meeting at a vertex can continue there. By our assumptions,
the same segment cannot be the edge of two triangles in the tiling.
Therefore, every uncut edge must continue at least at one of its endpoints.

We will use the following simple observation several times.

\medskip\noindent{\bf Lemma 5.} {\em If an edge $AB$ is subdivided and does not continue at $A$, then there is a unique triangle $ADE$ in the tiling such that $D$ is an interior
point of $AB$. The edge $AD$ is uncut. $\Box$}
\medskip

We call a triangle of the tiling {\it small} if all three of its sides are
uncut. We call it {\it large} if all three sides are subdivided. We call it
{\it improper} if it has an uncut side and also a side that does not
continue in either direction.

\medskip

 \begin{figure}[ht]
\centering
\includegraphics[height=3.5in]{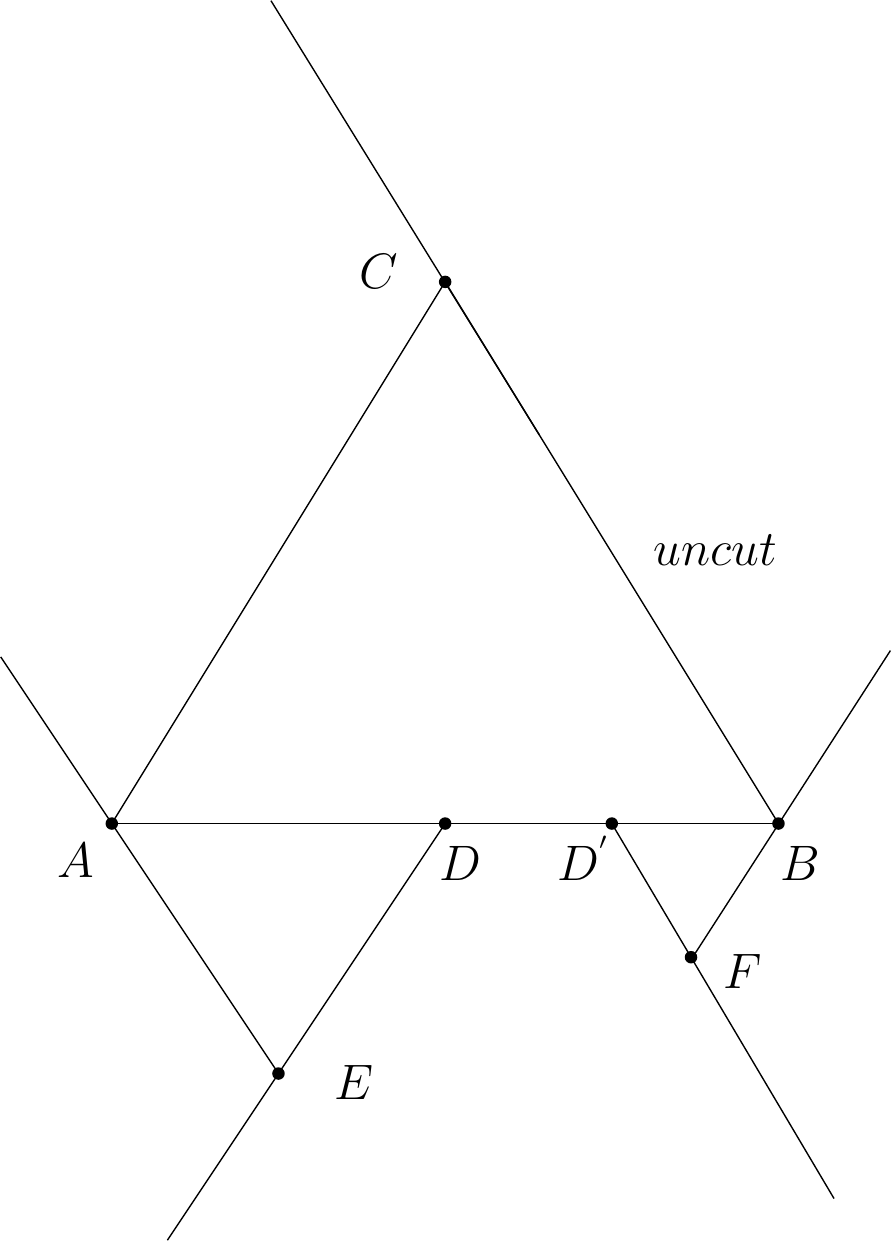}
\caption{Illustration to Lemma 6.}
\end{figure}

\medskip\noindent{\bf Lemma 6.} {\em For any improper triangle $T\in \cal T$, there are two other triangles, $U, V\in\cal T$ such that $U$ is improper
and the sum of the side lengths of $U$ and $V$ does not exceed the side length of $T$.}
\medskip

\noindent{\it Proof}: let $T=ABC$ and assume that the edge $BC$ is uncut and
$AB$ does not continue at either endpoint. As we noted, $AB$ must be
subdivided in this case and, by Lemma~5, we can find distinct
triangles $ADE$ and $D'BF$ such that $D$ and $D'$ are (not necessarily different) interior points of the
segment $AB$. See Figure~1. We clearly have $AD+D'B\le AB$ and we are done,
unless neither of these triangles is improper. Assume that they are not
improper. Their edges $AD$ and $D'B$ are
uncut, so all their edges have to continue in at least one direction. The edge
$AD$ continues at $D$, so $DE$ cannot continue there, so it must continue at
$E$. Therefore, the edge $AE$ must continue at $A$. Similarly, $FB$ must
continue at $B$. As $BC$ is uncut it must continue at one of its endpoints. As $FB$ continues at $B$, the edge $BC$ cannot continue there, so it
must continue at $C$. The edge $AC$ does not continue in
either direction, since $AE$ continues at $A$ and $BC$ continues at $C$.

As in the previous paragraph, we can find the triangles $AGH$ and $G'CI$ in the
tiling with $G,G'$ interior points of the segment $AC$ and we are done, unless
neither of them is improper. But, as we have seen earlier, in this case $IC$ would
continue at $C$, contradicting the fact $BC$ continues there. This contradiction
completes the proof of the lemma. $\Box$

\medskip\noindent{\bf Lemma 7.} {\em There is no improper triangle in the tiling.}
\medskip

\noindent{\it Proof}: Suppose $T_0$ is improper. By Lemma 6, for every $i\ge 1$ we can
recursively find triangles $T_i$ and $V_i\in\cal T$ such that $T_i$ is
improper and the sum of the side lengths of $T_i$ and $V_i$ does not
exceed the side length of $T_{i-1}$. Since the sum of the side lengths of
all triangles $V_i$ is at most the side length of $T_0$, this sum must be convergent, contradicting our assumption that the sizes of all triangles in $\cal T$ are separated from zero. The contradiction proves the
lemma. $\Box$

\medskip\noindent{\bf Lemma 8.} {\em If the edge $AB$ of a triangle $ABC$ in
$\cal T$ continues at $B$, then $BC$ is uncut.}
\medskip

\noindent{\it Proof}: As $AB$ continues at
$B$, $BC$ does not continue there, so if it is subdivided, then we
have a triangle $BDE$ such that $D$ is an interior point of $BC$ (Lemma~5). The edge
$BE$ does not continue at $B$ ($AB$ continues there), and $DE$ does not
continue at $D$ ($BD$ continues there). At most one of $BE$ and $DE$
can continue at $E$, so the triangle $ADE$ must have at least one edge that does
not continue in either direction. It also has an uncut edge, $BD$, so it is
improper. The contradiction with Lemma~7 proves Lemma~8. $\Box$

\medskip\noindent{\bf Lemma 9.} {\em Every triangle in $\cal T$ is either large or small. The
sides of the large triangles do not continue in either direction, and each of them
contains in its interior precisely one vertex. }
\medskip

\noindent{\it Proof}:
Assume first that the triangle $ABC$ has an edge that continues in at least
one direction, say $AB$ continues at $B$. By Lemma~8, $BC$ is uncut, so
it must also continue in one direction. It cannot continue at $B$, hence it
continues at $C$. By the same lemma, $AC$ must be uncut and it must
continue at $A$. Finally, using Lemma~8 again, we obtain that $AB$ is uncut and,
therefore, $ABC$ is a small triangle.

Assume next that none of the edges of $ABC$ continues in either
direction. Then all edges must be subdivided and $ABC$ is a large triangle. It
remains to prove that all three edges contain exactly one subdividing
vertex. Indeed, if $AB$ contains more than
one, then by an observation similar to Lemma~5, we have a triangle $DEF$ where
both $D$ and $E$ are interior vertices of the edge $AB$. In this case, $DF$ and
$EF$ do not continue at $D$ and $E$, respectively, ($DE$ continues at both
vertices), and at most
one of $DF$ and $EF$ can continue at $F$. So, $DEF$ has an edge that does not
continue in either direction. It also contains an uncut edge $DE$, so it is
improper. The contradiction with Lemma~7 completes the proof of
Lemma~9. $\Box$
\smallskip

 \begin{figure}[ht]
\centering
\includegraphics[height=3.5in]{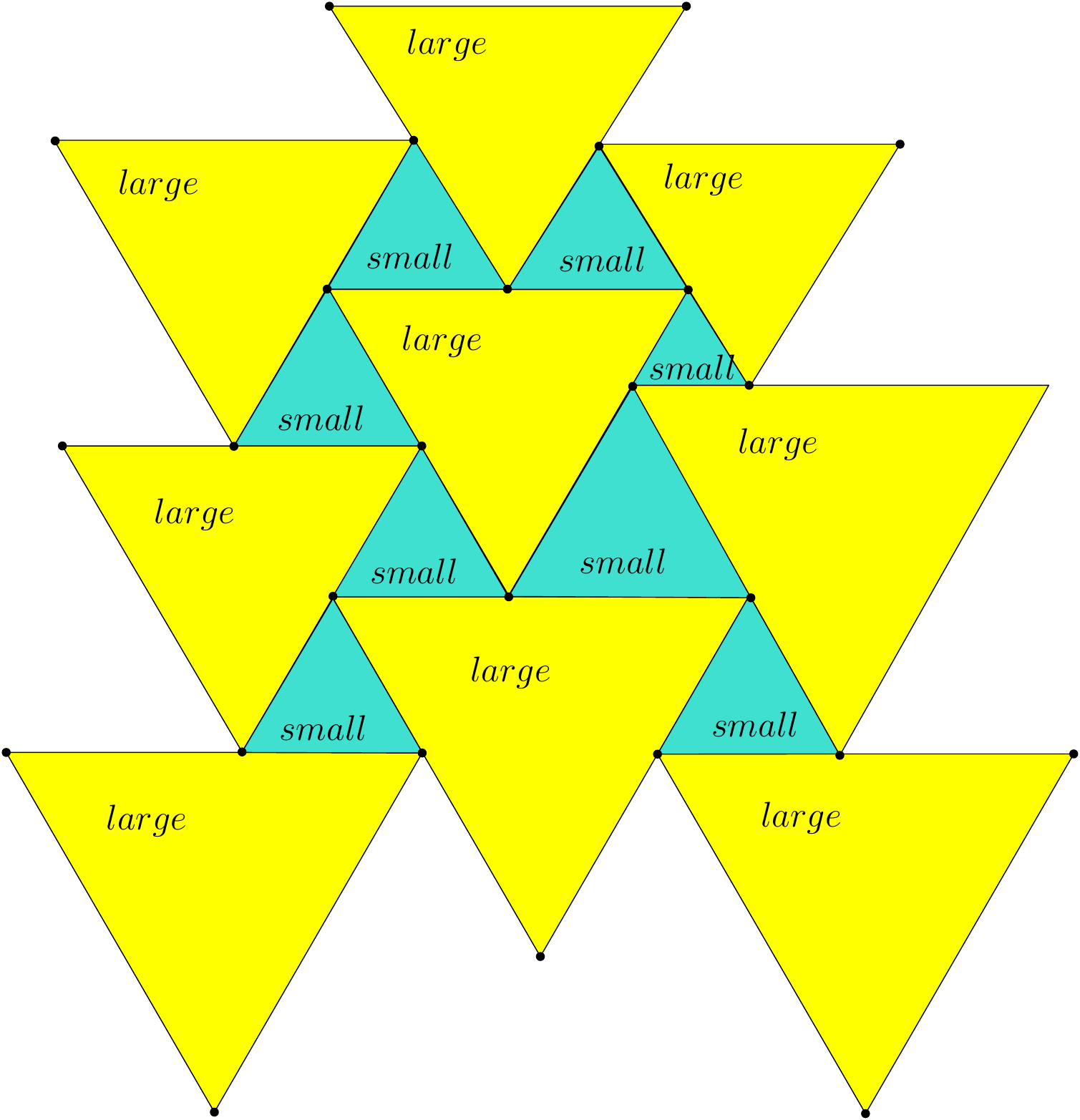}
\caption{The combinatorial structure of $\cal T$.}
\end{figure}

\smallskip

According to Lemma~9, the combinatorial structure of the
triangulation $\cal T$ is completely determined. A small part of the tiling is depicted in Figure~2. Each vertex of the
triangulation is the vertex of exactly one large triangle and exactly two
small triangles, and it is contained in exactly one further triangle: it is an interior vertex of an edge of a large triangle.

The following observation is immediate, using the fact that no triangles in the
tiling are improper.

\medskip\noindent{\bf Lemma 10.} {\em  Let $T=ABC$ be a large triangle in the tiling and $T_1$,
$T_ 2$ and $T_3$ be the large triangles containing $A$, $B$ and $C$,
respectively, in the interior of one of their edges. Then the side length of
$T$ is the average of the side lengths of $T_1$, $T_2$ and $T_3$. $\Box$}
\medskip

\noindent{\it Proof of Theorem 3}: Consider two large triangles $T$ and
$T'$. We claim that they must be congruent.

Assume for contradiction that the side length $a$ of $T$ is
smaller than the side length $a'$ of $T'$. Start a
random walk at $T_0=T$. After visiting the large triangle $T_i$, choose the next large triangle $T_{i+1}$ uniformly at random from among the three
triangles that contain a vertex of $T_i$ in the interior of one of its edges.
Stop the random walk at the first time when $T_i=T'$ holds.

According to a well known extension of P\'olya's recurrence theorem~\cite{Po21}, this random walk is recurrent; see~\cite{Ka97}, Theorem~8.2 on p.~138. Thus, choosing $N$ large enough, we can make sure that after $N$ steps it reaches $T'$ with probability at least
$a/a'$. We stop the process at step $N$,
even if $T'$ has not been reached.

The side lengths of $T_i$ form a martingale by Lemma~10. So the expected side
length at the end of the process is exactly the side length $a$ of $T_0=T$. At
the end, this side length is $a'$ with probability at least $a/a'$ and it
is always positive, so we have $a>(a/a')a'$. This contradiction proves that all large triangles must be of the same size.

In view of Lemma~8, every side of a large triangle is subdivided into two segments belonging to two small triangles. Since the structure of the tiling is uniquely determined, the small triangles must fall into at most two congruence classes, and the tiling is periodic (in fact, it is the union of three latticelike arrangements; see~\cite{Gr07}). $\Box$

\end{document}